\documentstyle[amscd]{amsart}
\input{bull-art}
\newtheorem{thm}{Theorem}
\newtheorem{cor}[thm]{Corollary}

\newtheorem{prop}[thm]{Proposition}

\ifx\undefined\theoremstyle
	\newtheorem{definition}[thm]{Definition}
	\newenvironment{defin}{\begin{definition}\rm}{%
\end{definition}}
	\newtheorem{conj}[thm]{Conjecture}
	\newtheorem{rem}[thm]{Remark}
	\newtheorem{rems}[thm]{Remarks}

\else
	\theoremstyle{definition}

	\newtheorem{defin}[thm]{Definition}

	\theoremstyle{remark}

	\newtheorem{rem}{Remark}

\fi

\ifx\undefined\pf
\newenvironment{pf}{\bigskip{\em Proof:\/}}{\qed\medskip}
\newenvironment{pf*}[1]{\bigskip{\em #1:\/}}{\qed\medskip}
\fi

\ifx\undefined\numberwithin
\def\numberwithin#1#2{\makeatletter\@ifundefined{c@#1}{%
\@nocnterrr}{%
  \@ifundefined{c@#2}{\@nocnterr}{%
  \@addtoreset{#1}{#2}%
  \toks@\expandafter\expandafter\expandafter{\csname 
the#1\endcsname}%
  \expandafter\xdef\csname the#1\endcsname
    {\expandafter\noexpand\csname the#2\endcsname
     .\the\toks@}}}\makeatother}\fi

\ifx\undefined\qed\newcommand{\qed}{\hfil\rule{4pt}{6pt}%
\bigskip}\fi
\newcommand{\Hom}{\operatorname{Hom}}  

\newcommand{\Map}{\mbox{\rm Map}}
\newcommand{\Aut}{\mathop{\mbox{\rm Aut}}}

\newcommand{\QZ}{{\Bbb Q}/{\Bbb Z}}

\newcommand{\Z}{\Bbb Z}

\newcommand{\Qp}{{\Bbb Q}_p}
\newcommand{\Zp}{{\Bbb Z}_p}
\newcommand{\zerowidth}[1]{\hbox to 0pt{\hss$\displaystyle 
#1$\hss}}
\newcommand{\LL}{[\mkern-2mu[}
\newcommand{\RR}{]\mkern-2mu]}

\ifx\undefined\eqref\newcommand{\eqref}[1]{\rm 
(\ref{#1})}\fi

\newcounter{thmItem}

\newenvironment{thmList}{\begin{list}%
{\rm \roman{thmItem})}{\usecounter{thmItem}
\setlength{\labelwidth}{2em}
\setlength{\itemindent}{2em}
\setlength{\leftmargin}{0pt}
\setlength{\listparindent}{0pt}
\setlength{\parsep}{0pt}
\setlength{\partopsep}{0pt}
\setlength{\itemsep}{\medskipamount}
\setlength{\topsep}{\medskipamount}
}}{\end{list}}

\newcounter{textItem}

\newcounter{condItem}

\theoremstyle{remark}

\newcommand{\Ga}{{\widehat{\Bbb G}_a}}
\newcommand{\Ext}{\operatorname{Ext}}

\newcommand{\gl}{\operatorname{GL}}
\newcommand{\liegl}{\operatorname{\mbox{$\frak g\frak l$}}}
\renewcommand{\P}{{\Bbb P}}
\newcommand{\pvk}{\P(V_K)}
\newcommand{\lk}{{\Cal L}_K}
\renewcommand{\l}{{\Cal L}}

\newcommand{\wseries}{W\LL u_1,\dots,u_{n-1}\RR}

\newcommand{\lie}{\operatorname{\mbox{$\Cal L$ie}}}
\newcommand{\M}{{\Cal M}}
\newcommand{\X}{{X}}

\newcommand{\lkn}{L_{K(n)}}
\newcommand{\morK}[1]{{\Cal K}_n({#1})}
\newcommand{\morM}[1]{{\Cal M}_n({#1})}
\newcommand{\hcts}{H_{\text{cts}}}
\newcommand{\Pic}{\operatorname{Pic}}
\newcommand{\K}{{\Cal K}}
\newcommand{\ox}{{\Cal O}_{\X}}
\newcommand{\oxk}{{\Cal O}_{\X_K}}
\newcommand{\C}[1]{{\Cal C}_{#1}}

\newcommand{\g}{{\frak g}}
\newcommand{\Der}{\operatorname{Der}}
\newcommand{\FF}{{\Cal F}}
\newcommand{\N}{{\operatorname{det}}}

\renewcommand{\Aut}{\operatorname{Aut}}

\begin{document}
\def\currentvolume{30}
\def\currentissue{1}
\def\currentyear{1994}
\def\currentmonth{January}
\def\copyrightyear{1994}
\def\currentpages{76-86}

\title[Lubin-Tate space and projective space]{The rigid 
analytic 
period mapping, Lubin-Tate space, and stable homotopy 
theory}
\ratitle
\author{M.~J.~Hopkins}
\address{Department of Mathematics \\ Massachusetts 
Institute of
Technology\\Cambridge, Massachusetts 02139}
\email{mjh@@math.mit.edu}
\author{B.~H.~Gross}
\address{Department of Mathematics \\ Harvard University\\ 
Cambridge,
Massachusetts 02138}
\email{gross@@math.harvard.edu}
\date{August 22, 1992}

\maketitle

\begin{abstract}
The geometry of the Lubin-Tate space of deformations of a 
formal group
is studied via an \'etale, rigid analytic map from the 
deformation
space to projective space.  This leads to a simple 
description of the
equivariant canonical bundle of the deformation space 
which, in turn, yields a
formula for the dualizing complex in stable homotopy theory.
\end{abstract}


\section*{Introduction}
Ever since Quillen~\cite{Qui:FGL,Ad:SHGH} discovered the 
relationship
between formal groups and complex cobordism, stable 
homotopy theory
and the theory of formal groups have been intimately 
connected.  Among
other things the height filtration of formal groups has 
led to the
chromatic filtration~\cite{Morava,MRW,Rav:Loc} which 
offers the best
global perspective on stable homotopy theory available.  
From the
point of view of homotopy theory this correspondence has 
been largely
an organizational principle.  It has always been easier to 
make
calculations with the algebraic apparatus familiar to 
topologists.
This is due to the fact that the geometry which comes up 
in studying
formal groups is the geometry of {\em affine} formal 
schemes.

The point of this paper is to study the Lubin-Tate 
deformation spaces
of 1-dimensional formal groups of finite height using a 
$p$-adic
analogue of the classical period mapping.  This is a rigid 
analytic,   
\'etale morphism, from the Lubin-Tate deformation space to 
projective
space, which is equivariant for the natural group action.  
With this
morphism the global geometry of projective space can be 
brought to
bear on the study of formal groups.  When applied to 
stable homotopy
theory, this leads to a formula for the analogue of the
Grothendieck-Serre dualizing complex.  

\section{Formal groups}

\subsection{The map to projective space}
Let $k$ be an algebraically closed field of characteristic 
$p>0$, and
let $F_0$ be a formal group of dimension 1 and finite 
height $n$
over~$k$.  Lubin and Tate~\cite{LT:Mod} studied the 
problem of deforming 
$F_0$ to a formal group $F$ over $R$, where $R$ is a 
complete, local
Noetherian ring with residue field $k$.  They defined 
deformations $F$
and $F'$ to be equivalent if there is an isomorphism
$
\phi:F\to F'
$
over $R$ which reduces to the identity morphism of $F_0$ 
and showed
that the functor which assigns to $R$ the equivalence 
classes of
deformations of $F_0$ to $R$ is representable by a smooth 
formal scheme $\X$
over the ring $W$ of Witt vectors of $k$.  Let $\Ga$ be 
the formal
additive group.  The essential cohomological
calculations in their argument are:
\begin{enumerate}
\item[(1)] $\Ext^2(F,\Ga)= 0$, which assures that 
deformations from $R/I$ to $R$
exist when $I^2=0$;
\item[(2)] $\Ext^1(F,\Ga)$ is a free $R$-module of rank 
$=n-1$, which determines the dimension
of the tangent space of the versal deformation;
\item[(3)] $\Hom(F,\Ga)=0$, which gives the existence of 
the moduli space.
\end{enumerate}

The space $\X$ is (noncanonically) isomorphic to the 
formal spectrum of
the power series ring $W\LL u_1,\dots,u_{n-1}\RR$ in $(n-1)$
variables  over $W$.  The group of automorphisms of the 
special fibre,
$
G=\Aut(F_0),
$
acts as formal automorphisms of
$\X$, as it acts as automorphisms of the functor 
represented by
$\X$~\cite{LT:Mod}.  An automorphism $T$ of $F_0$ deforms 
to an
isomorphism 
$
T:F_a\to F_b
$
of liftings over $R$ (corresponding to points $a,b\in 
\X(R)$)
if and only if $T\cdot a=b$.

Let $K$ be the quotient field of $W$, and let $\X_K$ be 
the rigid
analytic space over $K$ which is the ``generic fibre'' of 
the formal
scheme $\X$.  Then $\X_K$ is (noncanonically) isomorphic 
to the open
unit polydisc of dimension = $(n-1)$ over $K$, and $G$ 
acts as rigid
analytic automorphisms of $\X_K$.

The algebra 
$
D=\operatorname{End}(F_0)\otimes\Qp
$
is isomorphic to the division algebra of invariant 
$=\frac1n$ over
$\Qp$, and $\operatorname{End}(F_0)$ is the maximal 
$\Zp$-order in $D$.
Hence $G$ is isomorphic to the group of units of this 
order, which is
a maximal compact subgroup of the group $D^\times$.

Since the simple algebra
\[
D\otimes_{\Qp}K=M_n(K)
\]
is split, $G$ has a natural $n$-dimensional linear 
representation
$V_K$ over $K$.  It follows that $G$ acts by projective 
linear
transformation on the projective space $\pvk$ of 
hyperplanes in $V_K$.
This action extends to a projective linear action of the 
group $D^\times$.

Our main result is the following:

\begin{thm}\label{thm-main}
The crystalline period mapping \RM(to be defined below\RM) 
is an \'etale,
$G$-equivariant, rigid-analytic morphism
$
\Phi:\X_K\to \pvk.
$
Let $L$ be a Tate algebra \RM(= affinoid\RM) over $K$, and 
let $R\subseteq L$ be
the $W$-subalgebra of integral elements.  Let $F_a$ and 
$F_b$ be deformations
over $R$ of $F_0$ corresponding to the points 
\[
a,b\in\X_K(L)=\X(R).
\]
An isogeny $T$ of $F_0$ \RM(viewed as an element of 
$D^{\times}$\RM) deforms
to an isogeny
\[
f_T:F_a\to F_b
\]
if and only if
\[
T\cdot\Phi(a)=\Phi(b).
\]
\end{thm}

\begin{rem}
We have called $\Phi$ the crystalline period mapping to 
stress its
analogy with the period mapping from a simply connected 
complex family
of abelian varieties of complex dimension $g$ to the 
Grassmannian of
maximal isotropic subspaces in a complex symplectic space 
of dimension
$2g$ \cite{Grif:periods}.  Indeed both the crystalline and 
classical period
mappings are given by solutions to the Picard-Fuchs 
equation, which is the
differential equation given by the Gauss-Manin connection 
on the
primitive elements in the first deRham cohomology groups 
of the fibres.  In
the complex case the image lies in an open orbit (Siegel 
space) for
the real symplectic group by Riemann's positivity 
conditions.  In our
case the map $\Phi$ is surjective on points with values in 
the
completion of an algebraic closure of $K$.  The
fibres of $\Phi$ can be identified~\cite[\S23]{GH} with 
the cosets of
$\gl_n(\Zp)$ in the group
\[
\left\{\,g\in \gl_n(\Qp)\mid \det g\in\Zp^{\times}\right\}.
\]
\end{rem}

The proof of Theorem~\ref{thm-main} and its extension to 
formal $A$-modules
in the sense of Drinfeld~\cite{Drin:EM} are given in our 
paper~\cite{GH}.
There we provided explicit formulae for $\Phi$ using 
coordinates on
projective space;  here we will sketch a more abstract 
approach.  To
specify a $G$-morphism
\[
\Phi:\X_K\to\pvk,
\]
we must specify a $G$-equivariant, rigid analytic, line 
bundle $\lk$
on $\X_K$, as well as a homomorphism of $K[G]$-modules
\[
V_K\to H^0(\X_K,\lk)
\]
whose image has no base points.  Then $\Phi(x)$ is the 
hyperplane
$W_x$ in $V_K$, mapping to sections of $\lk$ which vanish 
at $x$.

Let $F$ be the universal deformation of $F_0$ over $\X$, 
and let
$e:\X\to F$ be the identity section of the formal morphism 
$\pi:F\to
X$.  Then
\[
\omega = \omega(F) = e^\ast\Omega^1_{F/\X}
\]
defines the invertible sheaf on $\X$ of invariant 
differentials on
$F$.  The sheaf
\[
\l=\omega^{-1}=\lie(F)
\]
is a $G$-equivariant line bundle on $\X$ whose fibres are 
the tangent
spaces of the corresponding deformations.  Let $\lk$ be 
the associated
rigid analytic line bundle on $\X_K$.  

Let $E$ be the universal additive extension of $F$ over 
$\X$;  then
$E$ is a formal group of dimension $=n$ which lies in an 
exact
sequence
\[
0\to N\to E\to F\to 0
\]
of formal groups over $\X$.  The formal group $N$ is 
isomorphic to the
additive group $\Ga\otimes\Ext(F,\Ga)^{\vee}$
 of dimension $=(n-1)$.  Passing to Lie
algebras gives an exact sequence
\[
0\to\lie(N)\to\lie(E)\to\lie(F)\to 0
\]
of $G$-equivariant vector bundles on $\X$.

The equivariant vector bundle
$
\M=\lie(E)
$
of rank $=n$ on $\X$ is the covariant Dieudonn\'e module 
of the formal
group $F$ (mod $p$).   As such $\M$ is a crystal over $\X$ 
in
the sense of Grothendieck~\cite{Groth:BT,MM}, it has an 
integrable connection
\[
\nabla:\M\to\M\otimes\Omega^1_{\X/W}
\]
together with a Frobenius structure.  A fundamental 
theorem of Dwork
on the radius of convergence of solutions to $p$-adic 
differential
equations (cf.\ \cite[Proposition~3.1]{Katz:Dwo}) shows 
that the vector
space
\[
V_K = H^0(\X_K,\M_K)^{\nabla}
\]
of rigid analytic, horizontal sections has dimension $n$ 
over $K$ and
that
\[
V_K\otimes\oxk\simeq\M_K.
\]
This vector space affords the natural $n$-dimensional 
representation of
$G$.  We stress that the horizontal sections in $V_K$ are 
rigid
analytic;  the $W$-module $H^0(\X,\M)^{\nabla}$ of formal 
horizontal
sections is zero once $n\ge 2$.

The surjection
\[
\M_K\to\lk\to 0
\]
of rigid analytic vector bundles gives a map of 
$K[G]$-modules
\[
V_K\to H^0(\X_K,\lk)
\]
whose image is base point-free.  This defines the 
crystalline period
mapping $\Phi$ in the theorem.  The fact that $\Phi$ is 
\'etale
follows from an explicit computation of the determinant of 
the
differential $d\Phi$~\cite[\S23]{GH}.

\begin{rem}
(1) The crystal $\M$ can also be identified with the 
primitive
elements in the first deRham cohomology group of
$F/X$~\cite{Katz:Bombay}.  From this point of view the
connection $\nabla$ is that of Gauss-Manin.

(2)
An analogue of the rigid analytic map $\Phi$ was 
introduced by
Katz~\cite{Katz:ICM} in the study of the moduli of 
ordinary elliptic
curves.  Katz called his map $L$ for logarithm.  The map 
$L$ was
studied in the supersingular case by Katz~\cite{Katz:Div} 
and
Fujiwara~\cite{Fujiwara}, and this is essentially 
equivalent to a
consideration of the map $\Phi$ in the case when $F_0$ has 
height
$n=2$.  In this case the inverse images under $\Phi$ of 
the two
points of $\P^1$ fixed by a maximal torus of $G$ are the 
moduli of
quasicanonical liftings of $F_0$, which were introduced
in~\cite{Gross:QC}.  In the general case quasicanonical 
liftings have
been studied by Jiu-Kang Yu~\cite{Yu}.

(3) 
The components of the map $\Phi$ were also investigated in 
joint
work of the first author with Ethan 
Devinatz~\cite{DevHop}.  That
account is presented in a language that might be more 
familiar to
topologists.
\end{rem}

\subsection{The $G$-action on $\X_K$ and the canonical 
line bundle}

The group $G$ is a (compact) $p$-adic Lie group.  Let $\g$ 
be its Lie
algebra over $\Qp$.  Then 
\[
\g\otimes K\simeq\liegl_n(K).
\]
If $\gamma\in\g$, then for $m\gg0$ the element 
$\exp(p^m\cdot\gamma)$
lies in $G$ and acts on $\X$.

\begin{prop}\label{prop-derivative}
Let $f$ be a rigid analytic function on $\X_K$, and let 
$\gamma$ be an
element of $\g$.  Then the limit
\[
D_\gamma(f) = \lim_{m\to\infty}\frac{\exp(p^m\gamma)\circ 
f-f}{p^m}
\]
exists in the Fr\'echet algebra $A_K$ of rigid analytic 
functions on
$\X_K$.  The map
$
f\mapsto D_\gamma(f\,)
$
is a derivation of $A_K$ over $K$, and the map
$
\gamma\mapsto D_\gamma
$
defines a representation of Lie algebras
$
\g\otimes K\to\Der_K(\X_K).
$
\end{prop}
Indeed, the corresponding facts are clear for the 
$G$-action on
$\pvk$, which is algebraic.  Since $\Phi$ is \'etale, 
there is no
obstruction to lifting the resulting vector fields to 
$\X_K$.
Informally speaking, the proposition shows that the action 
of $G$ on
$\X_K$ is ``differentiable'', which is not immediately 
apparent from
the definition of the $G$-action on $\X$.  In fact, if the 
ring of
formal functions $\ox$ is made into a normed algebra in 
the obvious way
(by choosing deformation parameters and giving them norm 1),
then the action map
\[
G\to\text{Bounded linear operators on $\ox$}
\]
is not even continuous.  It is not difficult to find 
explicit formulae
for the vector fields giving the differentiated
action~\cite[\S\S24 and 25]{GH}.

We can also use the map $\Phi$ to describe the canonical 
bundle of
$\X$ over $W$
\[
\Omega^{n-1}=\Omega^{n-1}_{\X/W}=\bigwedge^{n-1}\Omega^1_{%
\X/W}
\]
in the category of $G$-equivariant line bundles on $\X$.  

Let 
$\Theta=\Theta_{\X/W}$ be the tangent bundle of $\X$ over 
$W$.
The deformation theory of Kodaira and Spencer
(cf.\ [{7}, \S17; 13, Corollary~4.8] gives an isomorphism
of $G$-vector bundles
\[
\Theta \simeq\Hom(\lie(N),\lie(F)).
\]
Taking duals gives an isomorphism
$
\Omega^1\simeq\lie(N)\otimes\omega;
$
hence,
\[
\bigwedge^{n-1}\Omega^1\simeq\bigwedge^{n-1}\lie(N)\otimes%
\omega^{\otimes (n-1)}.
\]
Since the sequence
\[
0\to\lie(N)\to\lie(E)\to\lie(F)\to0
\]
of $G$-bundles is exact, we obtain a $G$-isomorphism
\[
\Omega^{n-1}\simeq\bigwedge^n\lie(E)\otimes\omega^{\otimes 
n}.
\]
It remains to identify the line bundle 
$\bigwedge^n\lie(E)$.  

If $\FF$ is an equivariant vector bundle on $\X$ and $k$ 
is an
integer, we let $\FF'=\FF[\N^k]$ be the equivariant bundle 
where the
action of $G$ on sections of $\FF$ is twisted by the
$k$th power of the reduced norm character
$\N:G\to\Zp^{\times}$:
\[
g'(f) = \N(g)^k\cdot g(f).
\]
By an analysis of the determinant of $d\Phi$, we 
show~\cite[\S22]{GH}
that $\Phi$ induces an isomorphism of equivariant line 
bundles on
$\X$:
\[
\bigwedge^n\lie(E)\overset{\sim}{\to}\ox[\N].
\]

Hence we obtain 

\begin{cor}\label{cor-canonical-bundle}
The canonical line bundle $\Omega^{n-1}$ of $\X$ is 
isomorphic, as a
$G$-equivariant line bundle, to 
$
\omega^{\otimes n}[\N].
$
\end{cor}

\begin{rem}
Since the map $\Phi$ is \'etale, this formula can also be 
deduced from
the (more elementary) corresponding formula for the 
$\gl_n$-equivariant bundle
$\Omega^{n-1}$ of $\pvk$.  Indeed,
$
\Omega^{n-1}_{\Bbb P}={\Cal O}_{\Bbb P}(-n)[\N]
$.
\end{rem}

\section{Duality in localized stable homotopy theory}
In this section we will work in the category of $p$-local 
spectra (in
the sense of stable homotopy theory).  Unfortunately, it 
does not seem
possible to preserve the standard notation of both 
algebraic geometry
and algebraic topology and avoid a conflict.  Throughout 
this section
the symbols $E$ and $F$ will denote spectra.  

The category of spectra is a triangulated category
and appears in the abstract to be much like the category 
of sheaves (or
complexes of sheaves) on a scheme.  For many purposes it 
is sufficient
to take this scheme to be a Riemann surface; but to 
understand the
more refined apparatus of stable homotopy theory, it is 
necessary to
take it to be a variety $S$ having a unique subvariety $S_n$
of each finite codimension $n$~\cite{Hop:GMHT,Morava}.

From the point of view of (complexes of) sheaves the 
sphere spectrum
corresponds to the sheaf of functions, the smash product 
of spectra
corresponds to the (derived) tensor product of sheaves, 
and stable
homotopy groups correspond to (hyper)cohomology groups.  A 
more
complete description of this analogy appears in Table 1.

\begin{table}\label{table-1}

\caption{Analogy between $p$-local spectra and quasicoherent
sheaves (or complexes of sheaves) on a variety $S$ having
a unique subvariety $S_n$ of each finite codimension $n$.  
The
map
$
i_n:U_n\to S
$
is the inclusion of $S\setminus S_n$} 

\begin{center}
\renewcommand{\arraystretch}{2}
\begin{tabular}{|p{1.85 in}||p{1.85 in}|}
\hline
\begin{center}\bf Stable Homotopy Theory\end{center} & 
\begin{center}\bf 
Sheaf Theory \end{center}\\
\hline
$p$-local spectrum & quasicoherent sheaf or complex of 
quasicoherent sheaves\\
``finite" $p$-local spectrum & finite complex of coherent 
sheaves \\
smash product  & tensor product \\
homotopy groups & hypercohomology \\
homotopy classes of maps & $\text{R\,Hom}$ \\
function spectra & sheaf $\text{R\,Hom}(A^*,B^*)$ \\
$p$-local sphere spectrum & $\Cal O_S$ \\
category $\C n$ & subcategory of coherent\newline 
sheaves~supported~on~$S_n$ \\
Morava $K(n)$ & total quotient field of \newline
$S_n\setminus\left(S_n\cap U_{n+1}\right)$ \\
functor $L_n$ & functor ${Ri_n}_\ast\circ{i_n}^\ast$\\
functor $\lkn$ & completion of $U_n$ along \newline 
$U_n\cap S_{n-1}$ \\
chromatic tower & Cousin complex \\
\hline
\end{tabular}
\end{center}
\end{table}

There are two kinds of duality in stable homotopy theory.  
The {\em
Spanier-Whitehead} dual of $F$, 
\[
DF=\Map[F,S^0],
\]
is the spectrum of maps from $F$ to the sphere spectrum.  
If $F$
is finite, then the homotopy type of $F$ is determined by 
a functorial
(in $E$) isomorphism
$
E^\ast F\approx E_{-\ast} DF.
$
The {\em Brown-Comenetz} dual $IF$ of $F$ represents the 
functor
\[
Y\mapsto\Hom\left(\pi_0 Y\wedge F,\QZ\right).
\]
If $I$ denotes the Brown-Comenetz dual of $S^0$, then 
there is a weak
equivalence
$
IF\approx\Map[F,I].
$

In the analogous situation of sheaves over a Riemann 
surface $S$, if
$F$ corresponds to a divisor $D$, then the 
Spanier-Whitehead dual of
$F$ corresponds to the divisor $-D$.  By the Serre duality 
theorem
the dual of $H^i(S,D)$ is $H^{1-i}(S,K-D)$, where $K$ is 
the canonical
sheaf.  It follows that the Brown-Comenetz dual $IF$ of $F$
corresponds to the divisor $K-D$ and that the dualizing 
complex $I$
corresponds to the divisor $K$ or, more precisely, to the 
complex of
sheaves consisting of the canonical sheaf in dimension 
$-1$ and zero
elsewhere.

The importance of this formula for the dualizing complex is
well known, and it is desirable to have an analogous 
description of
the dualizing complex $I$.  Since the category of finite 
spectra has
infinite Krull dimension~\cite{Hop:GMHT,HS}, it is 
necessary to
localize away from the primes of finite codimension.  This 
leads to
the chromatic tower.

\subsection{The chromatic tower}  In the situation of 
complexes of
sheaves over the filtered variety $S$ there is a standard 
procedure
for constructing a resolution of a sheaf by sheaves whose 
support lies
on one of the subvarieties $S_n\setminus S_{n+1}$.  This 
is called the
Cousin complex in~\cite{Hart:Res}.  In the context of 
stable homotopy
theory, it is known as the chromatic 
resolution~\cite{Rav:MU}.

Fix a rational prime~$p$, and let $L_n$ denote 
localization with
respect to the wedge 
$
K(0)\vee\dots\vee K(n)
$
of the first $n+1$ Morava $K$-theories~\cite{Rav:Loc}.  
There are natural
transformations
\begin{equation}\label{eq-chromatic-trans}
L_{n}\to L_{n-1}
\end{equation}
and compatible transformations
$
1\to L_n.
$
This results in the {\em chromatic tower}
\[
\begin{CD}@. \\
@. \vdots \\
@. L_1 S \\
@. @VVV\\
S\begin{picture}(0,0)
\put(5,12){\vector(1,1){29}}
\end{picture}
@>>> L_0 S\end{CD}
\]
When $F$ is the $p$-localization of a finite spectrum, the 
map
\[
\{\pi_k F\}\to\{\pi_kL_n F\}_n
\]
is a pro-isomorphism for each $k$~\cite{HR}.  This means 
that all of
the homotopy theory of finite spectra can be recovered 
from the
chromatic tower.

The difference between
$L_n$ and $L_{n-1}$ can be measured in two ways.  The 
fibre of the
transformation~\eqref{eq-chromatic-trans} is the functor 
$M_n$.  It is
known as the $n$th {\em monochromatic layer}.  The
difference is also measured by the functor $\lkn$, which 
is localization
with respect to the $n$th Morava $K$-theory.  It is believed
that the stable homotopy type of the $p$-localization of a 
finite
spectrum $F$ can in fact be recovered from the collective 
knowledge of
of the spectra $L_{K(n)}$.  There are natural equivalences
\begin{align*}
\lkn M_n F &\approx \lkn F,\\
M_n \lkn F&\approx M_n F,
\end{align*}
so the homotopy types of $\lkn F$ and $M_n F$ determine 
each other.

\subsection{The Morava correspondence}
As in \S1 fix a formal group of height $n\ge1$ over an
algebraically closed field $k$ of characteristic $p>0$, 
and let $G$ be
its automorphism group and $X$ the Lubin-Tate space of its
deformations.  

The homotopy types of the spectra $M_n F$ and $\lkn F$
are accessible through the Morava correspondence, which 
associates to
each spectrum~$F$, graded, equivariant, quasicoherent 
sheaves
$
\morK F\quad\text{and}\quad\morM F
$
over the Lubin-Tate space $\X$, and spectral sequences
\begin{align*}
\hcts^\ast(G;H^0\morK F) &\Rightarrow W\otimes\pi_*\lkn F, 
\\
\hcts^\ast(G;H^0\morM F) &\Rightarrow W\otimes\pi_*M_n F.
\end{align*}
These spectral sequences often collapse and always 
terminate at a
finite stage which depends only on $n$.  They converge to 
the graded
group associated to a finite filtration of their 
abutment~\cite{HR}.

When $F=S^0$ is the sphere spectrum, the
sheaf
$
\morK{S^0}
$
is the direct sum
$
\bigoplus_{n\in\Z}\l^{\otimes n}.
$
Its ring of global sections is
\[
\wseries[u,u^{-1}]\qquad \qquad 
\]
with $u\in H^0(\omega)$ a generator.  It is graded in such 
a way that
$|u|=-2$ and $|u_i|=0$.

Let 
\[
I_n=\Map[L_n S^0, I]
\]
be the Brown-Comenetz dual of $L_nS^0$.  There are maps 
$I_{n-1}\to
I_n$, and the direct limit
$
\varinjlim I_n
$
is the spectrum $I$.  The cofiber of $I_{n-1}\to I_n$ is 
$\lkn I_n$.
The spectrum $M_n I_n$ is the Brown-Comenetz dual of $\lkn 
S^0$.  
Knowledge of the homotopy type of $\lkn I_n$ completely 
describes
Brown-Comenetz duality in the category of $K(n)$-local 
spectra and, by
the above, ultimately gives a formula for the spectrum $I$.

\subsection{Invertible spectra}
\begin{defin}
A $K(n)$-local spectrum $F$ is {\em invertible} if there 
is a spectrum
$F'$ and a weak equivalence
\[
\lkn (F\wedge F')\approx\lkn S^0.
\]
\end{defin}

The invertible spectra form a group $\Pic_n$.

\begin{thm}\label{thm-invertible}
\RM{(i)} A $K(n)$-local spectrum $F$ is invertible if and 
only if the
associated sheaf $\morK F$ is the tensor product of an 
invertible
sheaf with $\morK{S^0}$.

\RM{(ii)}
If $p$ is large compared with $n$ $(2p-2\ge \max\{n^2,2n+
2\}$ will do\RM), then
the homotopy type of an invertible spectrum is determined 
by its
associated Morava module.
\end{thm}

\begin{thm} Let $\K$ be the total quotient field of the 
sheaf $\ox$.
The Morava modules associated to $\Sigma^{-n^2}I_n$ are
\begin{gather*}
\morM{\Sigma^{-n^2}I_n} = 
\K/\ox\otimes\Omega^{n-1}_{\X}\otimes\morK{S^0},\\
\morK{\Sigma^{-n^2-n}I_n} = 
\Omega^{n-1}_{\X}\otimes\morK{S^0}.
\end{gather*}
In particular, the spectrum $\lkn I_n$ is invertible.
\end{thm}

Together with Corollary~\ref{cor-canonical-bundle} this 
determines the
homotopy type of $\lkn I_n$ at large primes.  In certain 
cases 
(\S\ref{sec-2.4}) this
can be put into a more concrete form.

\subsection{Finite spectra and $v_n$ 
self-maps}\label{sec-2.4}

Let $\C0$ be the category of $p$-localizations of finite 
spectra, and
let
$
\C n\subseteq\C0
$
be the full subcategory of $K(n-1)$-acyclics.  These 
categories fit
into a sequence~\cite{HS}
\[
\dotsb\subset\C{n+1}\subset\C n\subset\dotsb\subset\C0.
\]
One of the main results of~\cite{HS} is that each $F\in\C 
n$ has, for $N\gg0$, 
an
``essentially unique'' $v_n$ self-map
\[
v:\Sigma^{2p^N(p^n-1)}F\to F
\]
satisfying
\[
K(m)_\ast v = \begin{cases}
0&\text{if $m\ne n$,} \\
\text{multiplication by $v_n^{p^N}$} &\text{otherwise}.
	      \end{cases}
\]
This gives for $N\gg0$ a canonical equivalence
\[
\lkn v:\lkn\Sigma^{2p^N(p^n-1)} F\overset{\sim}{\to}\lkn F.
\]
Because of this one can ``suspend'' $K(n)$-local spectra 
by any
element of 
\[
\varprojlim_N\,\Z/2p^N(p^n-1)\Z.
\]

\begin{cor}
Let $F$ be a spectrum in $\C n$. 
\begin{thmList}
\item[\RM{(i)}] If $p$ is odd and $n=1$, then
$
I_1F = \Sigma^2 DF.
$
\item[\RM{(ii)}] If $p$ is large with respect to $n$ \RM(as
in Theorem~\ref{thm-invertible}\RM) and $p1_F\sim\ast$, then
$
I_n F\approx\Sigma^\alpha\lkn  DF,
$
where
\[
\alpha = \lim_{N\to\infty}2p^{nN}\frac{p^n-1}{p-1}+n^2-n.
\]
\item[\RM{(iii)}] In particular, if $F$ admits a self-map 
$v$ satisfying
\[
K(n)_\ast v = v_n^{p^M},
\]
then there is a homotopy equivalence
\[
I_n F\approx \Sigma^{2p^{nM}(p^n-1)/(p-1)+n^2-n}DF.
\]
\end{thmList}
\end{cor}



\end{document}